\documentclass[10pt,twoside]{article}
\usepackage{amsmath,amstext,epsfig,graphicx,makeidx}
\usepackage{amscd,amsfonts,amssymb,amsthm}
\usepackage{t1enc}
 
 \newcommand{\g}[1][0]{\if#10
                         {\cal G}_n^\alpha
                         \else\if#11
                         {\cal G}_n^\alpha(z)
                         \else%
                         
1-(1-{#1})^{\alpha+2}\sum_{j=0}^{n-1}c_{\alpha,j}{#1}^j
                         \fi\fi}
 \newcommand{\apa}[1][0]{\if#10
                         \ifmmode A_\alpha^p\else$\!\!
A_\alpha^p$\fi
                         \else%
                         \ifmmode A_\alpha^{#1}\else$\!\!
A_\alpha^{#1}$\fi
                         \fi}
 \newcommand{\kka}[1][0]{\if#10
                         \ifmmode K_I^\alpha\else$K_I^\alpha$\fi
                         \else%
                         \ifmmode
K_{#1}^\alpha\else$K_{#1}^\alpha$\fi
                         \fi}
 \newcommand{\KA}[1][0]{\if#10
                         \ifmmode
K_A^\alpha(z,w)\else$K_A^\alpha(z,w)$\fi
                         \else%
                         \ifmmode
K_{#1}^\alpha(z,w)\else$K_{#1}^\alpha(z,w)$\fi
                         \fi}
 \newcommand{\kaf}[2][0]{\if#10
                         \ifmmode
K_A^\alpha(#2)\else$K_A^\alpha(#2)$\fi
                         \else%
                         \ifmmode
K_{#1}^\alpha(#2)\else$K_{#1}^\alpha(#2)$\fi
                         \fi}
 \newcommand{\fz}[1][0]{\if#10
                         (1-z\overline{w})^{\alpha+2}
                         \else%
                         (1-{#1})^{\alpha+2}
                         \fi}

 \newcommand{\cp}[1][0]{\if#10
                         \ifmmode \ce^+\else$\ce^+$\fi
                         \else\if#11
                         \ifmmode
\ce^+_{\epsilon}\else$\ce^+_{\epsilon}$\fi
                         \else%
                         \ifmmode \ce^+_{#1}\else$\ce^+_{#1}$\fi
                         \fi\fi}
 \theoremstyle{plain} 
 \newtheorem{thm}{THEOREM}[section] 
 \newtheorem{lemma}[thm]{LEMMA} 
 \newtheorem{prop}[thm]{PROPOSITION}
 
 \theoremstyle{definition} 
 
 \newcommand{\m}{{\mathbf K}}

\newcommand{\bk}{{\kappa\hspace{-.205cm}\kappa\hspace{-.21cm}\kappa\hspace{-.2105cm}\kappa}}

 \newcommand{\D}{{\mathbb D}} 
 \newcommand{\ii}{{\rm i}} 
 \newcommand{\A}{{\cal A}} 
 \newcommand{\re}{{\mathbb R}} 
 \newcommand{\ce}{{\mathbb C}}

 \newcommand{\T}{{\mathbb T}}
 \newcommand{\dis}{\displaystyle} 
  
 \newcommand{\real}{{\rm Re\:}} 
 \newcommand{\imag}{{\rm Im\:}}
 \newcommand{\cj}[1]{\overline{#1}}

 \newcounter{exampc}
 \newcommand{\theexamp}{\addtocounter{exampc}{1}
 \thesection.\arabic{exampc}}
 \newcommand{\example}[2][0]{\if#1%
                           0 {\footnotesize
                           \vspace{.3cm} \noindent {\bf Example
\theexamp}
                           {#2}\ \hfill $\diamond$  \vspace{.3cm}}
                           \else%
                           {\footnotesize \vspace{.3cm} \noindent
{\bf Example \theexamp (#1)}
                           {#2} \addvspace{0.3cm}}%
                           \fi}
 \newcommand{\ip}[3][0]{\if#10
             \left\langle {#2},{#3} \right\rangle
             \else\if#11
             \int_{\cp} {#2}(z)\cj{{#3}(z)}\omega(z)dA(z)
             \else
             \int_{#1} {#2}(z)\cj{{#3}(z)}\omega(z)dA(z)
             \fi\fi}

\numberwithin{equation}{section}
\title{Mean value surfaces with prescribed\\ curvature form}
\author{H\aa{}kan Hedenmalm \\ {\normalsize haakanh@math.kth.se}
\and Yolanda Perdomo G. \\ {\normalsize yolanda@maths.lth.se}}
\begin{document}
\maketitle

\noindent{\it Abstract:} The Gaussian curvature of a two-dimensional 
Riemannian manifold
is uniquely determined by the choice of the metric. The formulas
for computing the curvature in terms of components of the metric,
in isothermal coordinates, involve the Laplacian operator and
therefore, the problem of finding a Riemannian metric for a given
curvature form may be viewed as a potential theory problem. This
problem has, generally speaking, a multitude of solutions. To
specify the solution uniquely, we ask that the metric have the
mean value property for harmonic functions with respect to some
given point. This means that we assume that the surface is simply
connected and that it has a smooth boundary. In terms of the
so-called metric potential, we are looking for a unique smooth
solution to a nonlinear fourth order elliptic partial differential
equation with second order Cauchy data given on the boundary. We
find a simple condition on the curvature form which ensures that
there exists a smooth mean value surface solution. It reads: the
curvature form plus half the curvature form for the hyperbolic
plane (with the same coordinates) should be $\le0$. The same
analysis leads to results on the question of whether the canonical
divisors in weighted Bergman spaces over the unit disk have
extraneous zeros. Numerical work suggests that the above condition
on the curvature form is essentially sharp.

Our problem is in spirit analogous to the classical Minkowski
problem, where the sphere supplies the chart coordinates via the
Gauss map.
\medskip

\medskip

\medskip

\section{Introduction}

\noindent{\bf General abstract surfaces.} Consider a simply
connected $C^\infty$-smooth bordered two-dimensional Riemannian
manifold $\boldsymbol\Omega$. The boundary
$\partial{\boldsymbol\Omega}$ is then a $C^\infty$-smooth Jordan
curve. We model $\boldsymbol\Omega$ as the unit disk $\D$,
supplied with the $C^\infty$-smooth Riemannian metric $d{\mathbf
s}$:
\begin{equation}
{\boldsymbol\Omega}=(\D,d{\mathbf s}),\qquad d{\mathbf s}(z)^2
=a(z)\,dx^2+b(z)\,dy^2+2\,c(z)\,dxdy,
\label{eq-1}
\end{equation}
where we use the convention $z=x+iy$. The smoothness of the
metric means that $a,b,c$ are $C^\infty$-smooth real-valued
functions on the closed disk $\bar\D$, subject to the Riemannian
metric conditions
$$
0<a(z),b(z),\qquad c(z)^2<a(z)\,b(z).
$$
We are interested in the problem of reconstructing the
metric, if its associated curvature form is given. The {\sl area
form} is given by
$$
d{\boldsymbol\Sigma}(z)=\big(a(z)b(z)-c(z)^2\big)\,d\Sigma(z),
$$
where
$$
d\Sigma(z)=\frac{dxdy}{\pi},\qquad z=x+iy.
$$
The Gaussian {\sl curvature function} $\bk$ is a geometric
quantity, which is given by Brioschi's formula in terms of the 
coordinatization as
\begin{multline*}
\bk=\frac1{(a(z)b(z)-c(z)^2)^2}\Bigg\{
\\
\det\left(
\begin{array}{ccc} -\frac12 \partial_y^2 a(z) -\frac12
\partial_x^2 b(z)+\partial_{xy}^2 c(z) & \frac12 \partial_x a(z) &
\partial_x c(z)-\frac12 \partial_y a(z) \\
\partial_y c(z) -\frac12 \partial_x b(z) & a(z) & c(z) \\
\frac12\partial_y b(z) & c(z) & b(z)
\end{array}\right)
\\
-\det\left(\begin{array}{ccc}
0 & \frac12 \partial_y a(z) & \frac12 \partial_x b(z) \\
\frac12 \partial_y a(z) & a(z) & c(z) \\
\frac12 \partial_x b(z) & c(z) & b(z)
\end{array}\right)\Bigg\}
\end{multline*}
The {\sl curvature form} is given by the expression
$$
\m=\bk\:d{\boldsymbol\Sigma};
$$
it measures the distribution of the curvature in space. On the
curved surface $\boldsymbol\Omega$, there exists a counterpart
of the usual Laplacian in the plane, known as the {\sl
Laplace-Beltrami operator}, denoted by $\boldsymbol\Delta$. In
terms of the given coordinates, it can be expressed by
\begin{multline*}
\boldsymbol{\Delta} = \frac1{4\,\sqrt{a(z)b(z)-c(z)^2}}\bigg\{
\partial_x \bigg[\frac{b(z)}{\sqrt{a(z)b(z)-c(z)^2}}\,\partial_x
-\frac{c(z)}{\sqrt{a(z)b(z)-c(z)^2}}\,\partial_y\bigg]\\
+\partial_y\bigg[-\frac{c(z)}{\sqrt{a(z)b(z)-c(z)^2}}\,\partial_x
+\frac{a(z)}{\sqrt{a(z)b(z)-c(z)^2}}\,\partial_y\bigg]\bigg\}.
\end{multline*}
We say that a twice differentiable function $f$ on the curved
surface $\boldsymbol\Omega$ is {\sl harmonic} -- or
Laplace-Beltrami harmonic, if we want to emphasize that we use
the Laplacian induced by the metric -- provided that
${\boldsymbol\Delta}f(z)=0$ holds throughout
$\boldsymbol\Omega$. Now, suppose we have two curved surfaces
$\boldsymbol{\Omega}$ and ${\boldsymbol\Omega}'$, which are both
modelled by the unit disk $\D$:
$$
{\boldsymbol\Omega}=(\D,d{\mathbf s}),\qquad
{\boldsymbol\Omega}'= (\D,d{\mathbf s}').
$$
As the corresponding Laplace-Beltrami operators
$\boldsymbol\Delta$ and $\boldsymbol\Delta'$ are generally
different, we should expect them to give rise to different
collections of harmonic functions. If they give rise to identical
collections of harmonic functions, we say that the metrics
$d{\mathbf s}$ and $d{\mathbf s}'$ are {\sl isoharmic}. This
relation between two metrics of been isoharmic is an equivalence
relation and its equivalence classes are called {\sl isoharmic
classes of metrics}. Let $\mathfrak M$ denote the collection of
all $C^\infty$-smooth metrics on $\bar \D$, and let ${\mathfrak
M}_\alpha$ run through all the isoharmic classes as $\alpha$
passes through a suitably large index set. We then get the
disjoint decomposition
$$
{\mathfrak M} = \bigcup_\alpha {\mathfrak M}_\alpha,
$$
which provides a fibering of $\mathfrak M$. Of particular interest
is the fiber ${\mathfrak M}_0$ (assuming that the index set
contains the value 0) which contains the Euclidean metric as an
element. We claim that ${\mathfrak M}_0$ coincides with the
collection of so-called isothermal metrics. We recall the standard
terminology that $d{\mathbf s}$ is {\sl isothermal} if
\begin{equation}
d{\mathbf
s}(z)^2=\omega(z)\,|dz|^2=\omega(z)\,\big(dx^2+dy^2\big),
\label{eq-isoth}
\end{equation}
holds for some $C^\infty$-smooth function $\omega(z)$ which is
positive at each point of $\bar\D$. This means in terms of the
functions $a,b,c$ that
$$
a(z)=b(z)=\omega(z)\quad\text{and}\quad c(z)=0,
$$
hold throughout $\D$. The area form is then $d{\boldsymbol\Sigma}=
\omega d\Sigma$. The Laplace-Beltrami operator becomes a slight
variation of the ordinary Laplacian
$$
{\boldsymbol\Delta}=\frac1{\omega(z)}\,\Delta,
$$
where $\Delta$ denotes the normalized Laplacian
$$
\Delta = \Delta_z=\frac14\left(\frac{\partial^2}{\partial x^2} +
\frac{\partial^2}{\partial y^2}\right), \quad z=x+\ii y.
$$
This means that the Laplace-Beltrami harmonic functions for an
isothermal metric are just the ordinary harmonic functions. In
other words, the isothermal metrics are contained in ${\mathfrak
M}_0$. To see that all metrics in ${\mathfrak M}_0$ are
isothermal, we pick a metric $d{\mathbf s}$ in ${\mathfrak M}_0$
and note that the functions $f(z;z_0)=\real \big[(z-z_0)^2\big]$
and $g(z;z_0)=\imag \big[(z-z_0)^2\big]$ are ordinary harmonic
functions of $z$, and therefore also Laplace-Beltrami harmonic
functions. We obtain the equations ${\boldsymbol\Delta}_z
f(z;z_0)=0$ and ${\boldsymbol\Delta}_z g(z;z_0)=0$ and implement
the above formula for the Laplace-Beltrami operator
${\boldsymbol\Delta}$. We evaluate the Laplace-Beltrami equation
at the point $z_0$; as we vary the point $z_0$ in $\D$, the
assertion that the metric is isothermal follows.
\medskip

We would like to better understand the other isoharmic classes
${\mathfrak M}_\alpha$. To this end, we pick a metric $d{\mathbf
s}_1$ in ${\mathfrak M}_\alpha$, for a fixed index $\alpha$.
It is well-known in the theory of
quasi-conformal mappings that it is possible to find a
$C^\infty$-homeomorphism $F$ of $\bar \D$ which changes the
coordinate chart so that the metric $d{\mathbf s}_1$ becomes
isothermal \cite{ahlfors}. Under the same coordinate change $F$,
the other elements of ${\mathfrak M}_\alpha$ change as well.
However, the collection of Laplace-Beltrami harmonic functions,
being determined by the geometry of the abstract surface in
question, remains the same, except for the obvious composition
with the chart function $F$. The Laplace-Beltrami harmonic
functions for the metric $d{\mathbf s}_1$ are after the coordinate
change just the ordinary harmonic functions, and this then carries
over to all the other metrics in ${\mathfrak M}_\alpha$. In other
words, the $C^\infty$-homeomorphism $F$ effects an identification
of ${\mathfrak M}_\alpha$ with ${\mathfrak M}_0$.
\medskip

\noindent{\bf The optimization problem.} We are given a
$C^\infty$-smooth real-valued $2$-form $\boldsymbol\mu$ on the
closed unit disk $\bar\D$, and we are looking for a metric
$d{\mathbf s}$ on $\D$ of the type (\ref{eq-1}) which is smooth up
to the boundary and has $\boldsymbol\mu$ as curvature form
${\mathbf K}={\boldsymbol\kappa}\,d{\boldsymbol\Sigma}$. There are
plenty of such metrics. To reduce their number, we decide to
minimize the {\sl total area} of the associated surface
${\boldsymbol\Omega}= (\D,d{\mathbf s})$ under two additional
conditions. We require that (1) the area form at the origin is
essentially the area form of the plane, that is,
$$
a(0)b(0)-c(0)^2=1,
$$
and we also ask that (2) we only minimize over a fixed isoharmic
class ${\mathfrak M}_\alpha$ of metrics. This second requirement
means that we fix the collection of harmonic functions in the unit
disk while performing the minimization. A natural question here is
whether there exists a minimizing surface, and whether it has a
smooth boundary.

In order to solve this problem analytically, it is very helpful to
note that by using the $C^\infty$-homeomorphism $F$, mentioned in
the previous subsection, we may choose to work only with the class
${\mathfrak M}_0$ of isothermal metrics. This means that the
metrics are of the form
$$
d{\mathbf s}(z)^2
=\omega(z)\,|dz|^2=\omega(z)\,\big(dx^2+dy^2\big),
$$
for some $C^\infty$-smooth function $\omega(z)$ which is positive
at each point of $\bar\D$. The formula for the curvature
simplifies greatly,
$$
\bk(z)=-2{\boldsymbol\Delta}\log\omega(z)=
-\frac2{\omega(z)}\Delta\log\omega(z), \qquad z\in\D,
$$
and so does the formula for the curvature form:
\begin{equation}
{\mathbf
K}(z)=-2{\boldsymbol\Delta}\log\omega(z)\,d{\boldsymbol\Sigma}(z)
=-2\,\Delta\log\omega(z)\,d\Sigma(z),
\qquad z\in\D.
\label{eq-curv}
\end{equation}
Let us now see how this coordinatization simplifies the
formulation of our problem. We are given a $C^\infty$-smooth
real-valued function $\mu$ on the closed unit disk, which is the
density function for the $2$-form $\boldsymbol\mu$:
\begin{equation}
{\boldsymbol\mu}(z)=\mu(z)\,d\Sigma(z),\qquad z\in\D.
\label{eq-curvdata}
\end{equation}
We write the curvature form ${\mathbf K}$ as
$$
{\mathbf K}(z)=K(z)\,d\Sigma(z),\qquad
K(z)=-2\,\Delta\log\omega(z).
$$
The equation ${\mathbf K}={\boldsymbol\mu}$ then becomes
$$
-2\,\Delta\log\omega(z)=\mu(z),\qquad z\in\D.
$$
As we have already fixed the isoharmic class by setting the
Laplace-Beltrami harmonic functions equal to the ordinary
harmonic functions in the plane, all that we really need to fix
is the value of $\omega$ at the origin, while minimizing the
total area of the surface ${\boldsymbol\Omega}=(\D,d{\mathbf
s})$. However, {\sl it is equivalent to maximize the value
$\omega(0)$ while keeping the total area constant}, and we find
it convenient to fix the latter to equal $1$:
$$
|{\boldsymbol\Omega}|_{\boldsymbol\Sigma}=
|{\D}|_{\boldsymbol\Sigma}=
\int_\D d{\boldsymbol\Sigma}(z)=\int_\D \omega(z)\,d\Sigma(z)=1.
$$
\medskip

\noindent{\bf The optimization problem and the mean value
property.} We recall that the weight $\omega$ should solve the
following problem:
\begin{equation}
\left\{\begin{array}{l}
\mbox{maximize}\quad \omega(0),\quad\mbox{while}\\[.2cm]
\dis\Delta\log\omega(z)=-\frac12\,\mu(z),\qquad
z\in\D,\quad\mbox{and}
\\[.2cm]
\dis\int_\D \omega(z)\,d\Sigma(z)=1.
\end{array}
\right. \tag{OP}
\end{equation}

\noindent Here, $\mu$ is a $C^\infty$-smooth real-valued
function on $\bar\D$, and we ask of $\omega$ that it too should
be $C^\infty$-smooth (and positive) on $\bar\D$. It is not clear
that the optimization problem (OP) should have such a nice
solution in general. As a matter of fact, it is possible to
construct counterexamples that are fairly elementary.
Nevertheless, we shall investigate the properties that such an
extremal weight $\omega=\omega_0$ should enjoy. We compare the
extremal weight $\omega_0$ with nearby weights $\omega_t$, of
the form
$$
\omega_t(z)=e^{th(z)}\,\omega_0(z),
$$
where $h$ is a $C^\infty$-smooth real-valued function on
$\bar\D$ that is harmonic in the interior $\D$, with $h(0)=0$,
and $t$ is a real parameter. Then
$$
\Delta\log\omega_t(z)=\Delta\log\omega_0(z)=-\frac12\,\mu(z),
\qquad z\in\D,
$$
and $\omega_t(0)=\omega_0(0)$. The extremal property of
$\omega_0$ now forces the inequality
\begin{equation}
\int_\D\omega_0(z)\,d\Sigma(z)\le\int_\D\omega_t(z)\,d\Sigma(z)
=\int_\D e^{th(z)}\,\omega_0(z)\,d\Sigma(z)
\label{eq-2}
\end{equation}
to hold. By Taylor's formula,
$$
e^{th(z)}=1+th(z)+O(t^2),
$$
for $t$ close to $0$. As we plug this into equation
(\ref{eq-2}), we arrive at
\begin{equation}
\int_\D\omega_0(z)\,d\Sigma(z)\le\int_\D\omega_0(z)\,d\Sigma(z)
+t\int_\D h(z)\,\omega_0(z)\,d\Sigma(z)+O(t^2).
\label{eq-3}
\end{equation}
By varying $t$ from small positive to small negative values, we
realize that the only way for (\ref{eq-3}) to hold is if
\begin{equation*}
\int_\D h(z)\,\omega_0(z)\,d\Sigma(z)=0.
\end{equation*}
If we drop the requirement on $h$ that $h(0)=0$, and consider
the function $h(z)-h(0)$ instead in the above argument, we
obtain
\begin{equation*}
\int_\D h(z)\,\omega_0(z)\,d\Sigma(z)=h(0).
\end{equation*}
This is what we call the {\sl mean value property} of the weight
$\omega_0$. Note that by an approximation argument, the above
mean value property remains valid when we extend the collection
of $h$ to all harmonic functions in $\D$ that are integrable
with respect to area measure. We find that we are looking for a
(positive) weight $\omega_0$ that is $C^\infty$-smooth up to the
boundary, with
\begin{equation}
\Delta\log\omega_0(z)=-\frac12\,\mu(z),\qquad z\in\D,
\label{eq-4.0}
\end{equation}
and the mean value property
\begin{equation}
\int_\D h(z)\,\omega_0(z)\,d\Sigma(z)=h(0),\qquad h\in
{\mathcal H}^1(\D),
\label{eq-4}
\end{equation}
where ${\mathcal H}^1(\D)$ stands for the Banach space of all
complex-valued area integrable harmonic functions on $\D$.
Strictly speaking, there is more information contained in the
extremal property of $\omega_0$, but clearly, if we find a
unique $C^\infty$-smooth solution $\omega_0$ to (\ref{eq-4.0})
and (\ref{eq-4}), then it is definitely our prime candidate for
the solution to the optimization problem (OP).
\medskip

\noindent{\bf The relationship with Bergman kernel functions.}
We now discuss how to actually find the extremal weight
$\omega_0$; more precisely, we discuss the problem of solving
the equations (\ref{eq-4.0}) and (\ref{eq-4}). First, we note
that by elementary potential theory, one solution $\omega_1$ to
(\ref{eq-4.0}) is given by
\begin{equation}
\log\omega_1(z)=-\int_\D\log\left|\frac{z-w}{1-z\bar w}\right|\,
\mu(w)\,d\Sigma(w),\qquad z\in\D,
\label{eq-4.1}
\end{equation}
and it is well known that $\log\omega_1$ is real-valued and
$C^\infty$-smooth on $\bar\D$, because $\mu$ has these
properties. Any solution to (\ref{eq-4.0}), then, has the form
$$\log\omega_0(z)=\log\omega_1(z)+H(z),\qquad z\in\D,$$
where $H$ is real-valued and harmonic in $\D$; given the smoothness
assumptions on $\omega_0$, $H$ should be $C^\infty$-smooth on
$\bar\D$. We find a holomorphic function $F$ on $\D$, which is
zero-free and $C^\infty$-smooth up to the boundary, such that
$$
\log|F(z)|^2=H(z),\qquad z\in\D.
$$
By restricting (\ref{eq-4}) to holomorphic functions, we obtain
\begin{equation}
\int_D f(z)\,|F(z)|^2\omega_1(z)\,d\Sigma(z)=f(0),
\label{eq-6.111}
\end{equation}
for all $f$ in ${\mathcal A}^1(\D)$, the space of
area-integrable holomorphic functions on $\D$. Let
$g\in{\mathcal A}^1(\D)$ be arbitrary, except that $g(0)=0$;
then $f=g/F$ is in ${\mathcal A}^1(\D)$ as well, and we find
that (\ref{eq-6.111}) states that
\begin{equation}
\int_D g(z)\,\bar F(z)\,\omega_1(z)\,d\Sigma(z)=0.
\label{eq-6.112}
\end{equation}
We interpret this in terms of the Hilbert space ${\mathcal
A}^2(\D,\omega_1)$, consisting of the square area-integrable
holomorphic functions on $\D$, supplied with the weighted norm
$$
\|f\|_{\omega_1}=\left\{\int_\D |f(z)|^2\,
\omega_1(z)\,d\Sigma(z)\right\}^{1/2}.
$$
This space ${\mathcal A}^2(\D,\omega_1)$ is known as a {\sl
weighted Bergman space}. Equation (\ref{eq-6.112}) then states
that $F$ is perpendicular to all the functions
$$
\big\{g\in{\mathcal A}^2(\D,\omega_1):\,g(0)=0\big\},
$$
and this means that $F$ is of the form
$$
F(z)=C\,K_{\omega_1}(z,0),
$$
where $C$ is a complex constant, and $K_{\omega_1}(z,w)$ is the
{\sl weighted Bergman kernel with weight $\omega_1$}. We recall
that the weighted Bergman kernel is defined by
$$
K_{\omega_1}(z,w)=\sum_{n=1}^{+\infty}e_n(z)\,\bar e_n(w),\qquad
z,w\in\D,
$$
where the functions $e_1(z),\,e_2(z),\,e_3(z),\ldots$ run through
an orthonormal basis for ${\mathcal A}^2(\D,\omega_1)$ \cite[pp.
43-44]{BBergman}. Alternatively, $K_{\omega_1}(\cdot,w)$ is the
unique element of ${\mathcal A}^2(\D,\omega_1)$ which supplies the
point evaluation at $w\in\D$:
$$
f(w)=\int_\D f(z)\, \bar
K_{\omega_1}(z,w)\,\omega_1(z)\,d\Sigma(z),
$$
for $f\in{\mathcal A}^2(\D,\omega_1)$ \cite{Azj}. The constant $C$
is easily determined (at least in modulus) by applying
(\ref{eq-6.111}) with the choice $f=1$, which leads to
\begin{equation}
F(z)=K_{\omega_1}(0,0)^{-1/2}\,K_{\omega_1}(z,0),\qquad z\in\D.
\label{eq-6.113}
\end{equation}
Returning back to the extremal weight $\omega_0$, we find that
it is of the form
\begin{equation}
\omega_0(z)=\frac{\big|K_{\omega_1}(z,0)\big|^2}{K_{\omega_1}(0,0)}
\,
\omega_1(z),\qquad z\in\D.
\label{eq-6.114}
\end{equation}
By the elliptic regularity theory for PDEs \cite{Niren}, the
function $K_{\omega_1}(\cdot,0)$ is $C^\infty$-smooth on $\bar\D$.
We realize, then, that {\sl a necessary condition for the
existence of a $C^\infty$-smooth positive weight $\omega_0$ on
$\bar\D$ that solves our problem (OP) is that
\begin{equation}
K_{\omega_1}(z,0)\neq0,\qquad z\in\bar\D.
\label{eq-6.115}
\end{equation}
It turns out that it is also sufficient, and that the extremal
solution is then given by (\ref{eq-6.114})}. For, if $\omega$ is
another weight that solves (\ref{eq-4.0}), and is
$C^\infty$-smooth and positive on $\bar\D$, then it is of the
form
$$
\omega(z)=|F(z)|^2\,\omega_1(z),\qquad z\in\D,
$$
where $F$ is $C^\infty$-smooth on $\bar\D$, analytic in $\D$,
and zero-free in $\bar\D$. It is well-known that the function
$$
F(z)=K_{\omega_1}(0,0)^{-1/2}\,K_{\omega_1}(z,0),\qquad z\in\D,
$$
is the unique (up to multiplication by unimodular constants)
solution to the extremal problem to maximize $|F(0)|$, given
that $F$ is holomorphic in $\D$ and
$$
\int_\D |F(z)|^2\,\omega_1(z)\,d\Sigma(z)=1,
$$
which means that $\omega(0)\le\omega_0(0)$ for all competitors
$\omega$ with
$$
\int_\D \omega(z)\,d\Sigma(z)=1.
$$
This shows that $\omega_0$, as given by (\ref{eq-6.114}), is
indeed the solution to the extremal problem, provided that
condition (\ref{eq-6.115}) is fulfilled.
\medskip

\noindent{\bf A toy example.} We should study a simple example, to
develop some intuition. We consider the degenerate data
$$
\mu(z)=-\theta\,\delta_\lambda(z),
$$
where $\delta_\lambda(z)$ stands for the Dirac delta function
concentrated at the point $\lambda\in\D$, and $\theta$ is a real
parameter. This does not fulfill our smoothness requirement at
the point $\lambda$, so we should think of it as a limit case of
smooth functions $\mu_n$. The weight function $\omega_1$
supplied by equation (\ref{eq-4.1}) is explicitly given by
$$
\omega_1(z)=\left|\frac{z-\lambda}{1-\bar\lambda
z}\right|^\theta,\qquad z\in\D,
$$
which is a reasonable weight on $\D$ provided that
$-2<\theta<+\infty$; outside this interval, the weight fails to
be area-summable near the point $\lambda$. We need to find the
weighted Bergman kernel function $K_{\omega_1}$. To this end, we
note first that if $\phi$ is a M\oe{}bius automorphism of $\D$,
then we have the following relationship between the kernel
functions for the weights $\omega$ and $\omega\circ\phi$:
$$
K_{\omega\circ\phi}(z,w)=
\phi'(z)\,\bar\phi'(w)\,K_\omega\big(\phi(z),\phi(w)\big),
\qquad z,w\in\D.
$$
Let $\omega_2$ stand for the radial weight
$$
\omega_2(z)=|z|^\theta,\qquad z\in\D;
$$
then the weighted Bergman space ${\mathcal A}^2(\D,\omega_2)$
has the reproducing kernel function
$$
K_{\omega_2}(z,w)=\frac1{(1-z\bar
w)^2}+\frac{\theta}2\,\frac1{1-z\bar w}, \qquad z,w\in\D.
$$
As $\omega_1=\omega_2\circ\phi$, where $\phi$ is the involutive
M\oe{}bius automorphism
$$
\phi(z)=\frac{\lambda-z}{1-\bar \lambda z},\qquad z\in\D,
$$
we find that
$$
K_{\omega_1}(z,w)=\frac1{(1-z\bar w)^2}+\frac{\theta}2\,
\frac{1-|\lambda|^2}{(1-\bar\lambda z)(1-\lambda\bar w)(1-z\bar
w)}, \qquad z,w\in\D.
$$
Plugging in $w=0$, we obtain
$$
K_{\omega_1}(z,0)=1+\frac{\theta}2\,
\frac{1-|\lambda|^2}{1-\bar\lambda z}, \qquad z\in\D.
$$
This expression has a zero in $\D$ precisely when (recall that
$-2<\theta<+\infty$ is assumed)
$$
-2<\theta<-\frac2{1+|\lambda|}.
$$
In particular, if $-2<\theta<-1$, we may choose $\lambda$ close to
the unit circle $\T$, to make sure that $K_{\omega_1}(z,0)$ has a
zero in $\D$, whereas if $-1\le\theta<+\infty$, no such zero can
be found, no matter how cleverly we try to pick $\lambda\in\D$.
This means that the function $\omega_0$ supplied by relation
(\ref{eq-6.114}) is $C^\infty$-smooth on
$\bar\D\setminus\{\lambda\}$ for $-1\le\theta<+\infty$, which
constitutes the solution to our extremal problem. However, for
$-2<\theta<-1$, we see that it is possible to pick $\lambda\in\D$
so that no smooth solution to the extremal problem will exist.

This calculation suggests a pattern: {\sl for hyperbolic metrics
(this means that the Gaussian curvature is negative everywhere),
we have a smooth solution $\omega_0$ to the extremal problem,
whereas when the metric becomes elliptic (positive Gaussian
curvature), we tend to get in trouble}.
\medskip

\noindent{\bf Statement of the main results.} We are going to
compare the curvature form with that of the hyperbolic plane,
and use this as a criterion for how strongly curved our metric
is. The Poincar\'e metric is
$$
d{\mathbf s}_{\mathbb H}(z)=\frac{2\,|dz|}{1-|z|^2},\qquad
z\in\D,
$$
and the associated area form is
$$
d{\boldsymbol\Sigma}_{\mathbb H}(z)=
\frac{4\,d\Sigma(z)}{(1-|z|^2)^2}.
$$
The curvature form for the Poincar\'e metric (which supplies the
standard model for the hyperbolic plane $\mathbb H$) is
$$
{\mathbf K}_{\mathbb H}(z)=-\frac{4\,d\Sigma(z)}{(1-|z|^2)^2};
$$
we are to compare the curvature form data
$\boldsymbol\mu(z)=\mu(z)\,d\Sigma(z)$ with ${\mathbf
K}_{\mathbb H}(z)$.

\begin{thm}
Let $\mu$ be a $C^\infty$-smooth real-valued function on
$\bar\D$, and suppose that the associated $2$-form
$\boldsymbol\mu(z)=\mu(z)\,d\Sigma(z)$ has
$$
\boldsymbol\mu(z)+\frac12\, {\mathbf K}_{\mathbb H}(z)\le0,\qquad z\in\D.
$$
Then the optimization problem (OP) has a unique
$C^\infty$-smooth positive solution $\omega_0$ on $\bar\D$,
which is given by relation (\ref{eq-6.114}), where $\omega_1$ is
as in (\ref{eq-4.1}). In addition,
$$
K_{\omega_1}(z,w)\neq0,\qquad
(z,w)\in\big(\bar\D\times\D\big)\cup \big(\D\times\bar\D\big).
$$
\label{thm-main}
\end{thm}

\noindent This is in line with the intuition we arrived at from
our ``toy example''. Note that the assumption of the theorem is
considerably weaker than requiring negative data $\boldsymbol\mu$.
As for the necessity of the condition of Theorem \ref{thm-main},
we have obtained the following.

\begin{thm}
Fix $\alpha$, $\alpha_0<\alpha<+\infty$, where
$\alpha_0\approx1.04$. Consider a $C^\infty$-smooth real-valued
function $\mu$ on $\bar\D$, and suppose that the associated
$2$-form $\boldsymbol\mu(z)=\mu(z)\,d\Sigma(z)$ has
$$
\boldsymbol\mu(z)+\frac\alpha2\, {\mathbf K}_{\mathbb H}(z)\le0,\qquad
z\in\D.
$$
Then there exists a choice of $\mu$ such that the optimization
problem (OP) fails to have a $C^\infty$-smooth positive solution
$\omega_0$ on $\bar\D$. \label{thm-nec}
\end{thm}

We conjecture that Theorem \ref{thm-nec} will remain true with
$\alpha_0=1$, making the statement of Theorem \ref{thm-main}
essentially sharp.

Work related to the problems considered here can be found in the papers
\cite{DKS1994}, \cite{Hed1}, \cite{HHopen}, \cite{DKS1996},
\cite{HHfact}, \cite{SS1}, \cite{SS2}, \cite{H&Z}, \cite{Weir}.

\medskip

\noindent{\bf The metric potential.} The Green function for the
Laplacian $\Delta$ is
$$
G(z,w)=\log\left|\frac{z-w}{1-z\bar w}\right|^2,\qquad
z,w\in\D,\quad z\neq w.
$$
The {\sl metric potential} associated with the isothermal metric
(\ref{eq-isoth}) is the function
$$
{\boldsymbol\Phi}(z)=\int_\D G(z,w)\,\omega(w)\,d\Sigma(w),
\qquad z\in\D;
$$
see, for instance, \cite{Jost}. It solves the boundary value problem 
($\T$ is the unit circle)
\begin{equation}
 \left\{\begin{array}{l}
{\boldsymbol\Delta}{\boldsymbol\Phi}(z)=1,\qquad z\in\D,
  \\[0.1cm]
{\boldsymbol\Phi}(z)=0,\qquad z\in\T.
 \end{array}\right.
 \label{mpot}
\end{equation}
It is set in boldface because it expresses a quantity that is 
independent of the choice of coordinates. We wish
to describe the potential equation (\ref{eq-4.0}) and the mean
value property (\ref{eq-4}), which the smooth solution to the
optimization problem (OP) should satisfy, in terms of the metric
potential. If we let ${\boldsymbol\Phi}_0$ stand for the metric
potential associated with a weight $\omega_0$ with
(\ref{eq-4.0}) and (\ref{eq-4}), we find that
${\boldsymbol\Phi}_0$ solves
\begin{equation}
 \left\{\begin{array}{l}
\dis\Delta\log\Delta{\boldsymbol\Phi}_0(z)=-\frac12\,\mu(z),\qquad
z\in\D,
  \\[0.3cm]
\dis{\boldsymbol\Phi}_0(z)=0,\qquad z\in\T,\\[0.2cm]
\dis\frac{\partial}{\partial n(z)}\,{\boldsymbol\Phi}_0(z)=2,\qquad
z\in\T,
 \end{array}\right.
 \label{eq-mpot}
\end{equation}
where $\partial/\partial n(z)$ denotes the normal derivative,
taken in the exterior direction. We are of course only looking for
subharmonic solutions ${\boldsymbol\Phi}_0$, which means that the
expression $\log\Delta{\boldsymbol\Phi}_0$ is more or less
well-defined (at least if we have some additional smoothness, and
the subharmonicity is ``strong''). The normal derivative
condition in (\ref{eq-mpot}) cleverly encodes the mean value
property (\ref{eq-4}), as is seen easily from an application of
Green's formula. {\sl We realize that Theorem \ref{thm-main} can
be interpreted as an assertion claiming the existence of a unique
smooth solution to the non-linear elliptic boundary value problem
(\ref{eq-mpot}) under appropriate conditions on $\mu$}. By moving
the point at the origin around by applying a M\oe{}bius
transformation that fixes the unit disk, we find that the method
can also treat the case when the data for the normal derivative is
replaced by
$$
\frac{\partial}{\partial n(z)}\,{\boldsymbol\Phi}_0(z)=
2\,\frac{1-|\lambda|^2}{|\lambda- z|^2},
$$
where $\lambda$ is any point of $\D$. It would be interesting to
have an analysis of the equation (\ref{eq-mpot}), where the
boundary data are considered from a wider class of functions.


\section{Preliminaries}

\noindent{\bf The standard weighted Bergman spaces.} For
$-1<\alpha<+\infty$, let $\A_\alpha^2(\D)=\A^2(\D,\omega_\alpha)$
denote the weighted Bergman space for the weight
$$
\omega_\alpha(z)=(\alpha+1)(1-|z|^2)^\alpha,\qquad z\in\D.
$$
The norm in $\A_\alpha^2(\D)$ is written
$$
\|f\|_{\alpha}=\left\{\int_\D|f(z)|^2\,\omega_\alpha(z)\,d\Sigma(z)
\right\}^{1/2}.
$$
Take a point $\lambda\in\D\setminus\{0\}$, and let
$\varphi_\lambda(z)= \varphi_{\lambda,\alpha}(z)$ be the
so-called {\sl extremal function} for the problem
$$
\sup\Big\{\real f(0):\: f(\lambda)=0,\; \| f\|_\alpha\leq1\Big\},
$$
which is unique, by elementary Hilbert space theory. If we let
$K_\lambda=K_{\lambda}^{\alpha}$ denote the reproducing kernel for the
closed subspace
$$
\Big\{f\in\A_\alpha^2(\D):\,f(\lambda)=0\Big\},
$$
then the extremal function can be written
$$
\varphi_{\lambda}(z)=K_{\lambda}(0,0)^{-1/2}\,
K_{\lambda}(z,0),\qquad z\in\D.
$$
When the function $\varphi_{\lambda}$ serves as a good divisor of
the zero at $\lambda$ in the space $\A_\alpha^2(\D)$, it is called
the {\it canonical divisor} of $\lambda$ as it is called in
\cite{DKS1994}, or the {\it contractive zero divisor}, provided
that division by it defines a norm contractive operation as it is
called in \cite{Hed1}. It is known \cite[pp. 58]{HHbook} (see
Section 4 for details) that the reproducing kernel $K=K_\alpha$
for the space $\A_\alpha^2(\D)$ has the form
$$
K(z,w)=\frac{1}{(1-z\bar w)^{\alpha+2}},\qquad z,w\in\D,
$$
and that the kernel $K_{\lambda}$ is derived from $K$ via the
identity \cite{HHbook}
\begin{multline*}
K_{\lambda}(z,w)=
K(z,w)-\frac{K(z,\lambda)K(\lambda,w)}{
K(\lambda,\lambda)}
\\
=\frac{1}{(1-z\bar
w)^{\alpha+2}}-\frac{(1-|\lambda|^2)^{\alpha+2}}
{(1-z\bar\lambda)^{\alpha+2}(1-\lambda\bar w)^{\alpha+2}},
\end{multline*}
(see Section 4 for the details). It follows that the extremal
function is given by
\begin{equation}
\varphi_{\lambda}(z)=\big(1-
(1-|\lambda|^2)^{\alpha+2}\big)^{-1/2} \left\{1-
\frac{(1-|\lambda|^2)^{\alpha+2}}{(1-z\bar\lambda)^{\alpha+2}}
\right\}, \qquad z\in\D. %
\label{extremal function}
\end{equation}

The following property of the extremal function $\varphi_\lambda=
\varphi_{\lambda,\alpha}$ is fundamental.

\begin{lemma}\label{lemma 2.1}
Fix $\alpha$ in the interval $-1<\alpha<+\infty$. Then, for each
bounded
harmonic function $h$ on $\D$, we have that
$$
\int_\D
h(z)\,|\varphi_\lambda(z)|^2\,(\alpha+1)\,\big(1-|z|^2\big)^\alpha
\,d\Sigma(z) = h(0).
$$
\end{lemma}
\begin{proof}
First, suppose $h$ is an analytic polynomial.
Then, if we recall the definition of $\varphi_\lambda$ in terms of
$K_\lambda$,
and use the reproducing property of the kernel $K_\lambda$, we have
that
\begin{multline}
\int_\D h(z)\,|\varphi_\lambda(z)|^2(\alpha+1)\big(1-|z|^2\big)^\alpha
\,d\Sigma(z)
\\
= \int_\D h(z)\frac{K_{\lambda}(z,0)}{K_{\lambda}(0,0)}
\bar K_\lambda(z,0)\,(\alpha+1)\big(1-|z|^2)^\alpha d\Sigma(z) =h(0).
\label{eq-l2.1}
\end{multline}
Taking complex conjugates in (\ref{eq-l2.1}), we realize that the
desired equality holds for all harmonic polynomials (defined to be
sums of analytic and antianalytic polynomials). A simple
approximation argument finishes the proof.
\end{proof}

\medskip

\noindent{\bf A weighted biharmonic Green function.}
The biharmonic Green function is the function
$\Gamma$ on $\D\times\D$ that solves the boundary value problem
\begin{equation}
\left\{\begin{array}{l}
\dis\Delta_z^2\Gamma(z,w)=\delta_w(z),\qquad z\in\D,
  \\[0.3cm]
\dis\Gamma(z,w)=0,\qquad z\in\T,\\[0.2cm]
\dis\frac{\partial}{\partial n(z)}\,\Gamma(z,w)=0,\qquad z\in\T,
 \end{array}\right.
 \label{eq-biharm}
\end{equation}
where $w\in\D$, and $\delta_w$ stands for the unit point mass at
the point $w$. We will think of locally summable functions $f$ on some
domain $\Omega$ of the complex plane as distributions on $\Omega$ via
the
linear duality
$$\langle f,\phi\rangle=\int_\Omega f(z)\,\phi(z)\,d\Sigma(z),$$
for compactly supported test functions $\phi$ on $\Omega$.
Given this normalization, the biharmonic Green function is given
explicitly by the formula
$$\Gamma(z,w)=|z-w|^2\log\left|\frac{z-w}{1-z\bar w}\right|^2+
(1-|z|^2)(1-|w|^2),\qquad z,w\in\D.$$ We shall also need the
Green function for the weighted biharmonic operator
$\Delta(1-|z|^2)^{-1}\Delta$, denoted $\Gamma_1$, which, by
definition, solves, for fixed $w\in\D$,
\begin{equation}
 \left\{\begin{array}{l}
\dis\Delta(1-|z|^2)^{-1}\Delta\Gamma_1(z,w)=\delta_w(z),\qquad
z\in\D,
  \\[0.3cm]
\dis\Gamma_1(z,w)=0,\qquad z\in\T,\\[0.2cm]
\dis\frac{\partial}{\partial n(z)}\,\Gamma_1(z,w)=0,\qquad z\in\T,
 \end{array}\right.
 \label{eq-wbiharm}
\end{equation}
Although the differential operator is singular at the boundary,
the above boundary value problem has a unique solution (we may use
Green's theorem to interpret the boundary data in terms of
integral conditions for $\Delta_z\Gamma_1(z,w)$, which are
uniquely solvable).

The function $\Gamma_1$ was calculated explicitly in
\cite{HHopen}:
\begin{lemma}
We have that
\begin{multline*}
\Gamma_1(z,w)=\bigg\{|z-w|^2-\frac14\big|z^2-w^2\big|^2\bigg\}
\log\left|\frac{z-w}{1-z\bar w}\right|^2\\
+\frac18\,(1-|z|^2)(1-|w|^2)
\Bigg\{7-|z|^2-|w|^2-|zw|^2-4\real z\bar w \\
-2\,(1-|z|^2)(1-|w|^2)\frac{1-|zw|^2}{|1-z\bar w|^2}\Bigg\},
\end{multline*}
for $z,w\in\D$. Moreover, it follows that, for $z,w\in\D$,
\begin{multline*}
\frac18\,\frac{(1-|z|^2)^3(1-|w|^2)^3}{|1-z\bar w|^2}\le
\Gamma_1(z,w)\\
\le\frac18\,\frac{(1-|z|^2)^3(1-|w|^2)^3}{|1-z\bar w|^4}
\Big\{|1-z\bar w|^2+4-|z+w|^2\Big\}.
\end{multline*}
\end{lemma}

\begin{proof}
Let $\widetilde\Gamma_1$ stand for the function defined by the above
expression; we want to show that $\widetilde\Gamma_1=\Gamma_1$. To this
end, we show that it solves the boundary value problem which determines
$\Gamma_1$. We first note that
$$
\Delta_z\widetilde\Gamma_1(z,w)=(1-|z|^2)\big[G(z,w)+H_1(z,w)\big],
\qquad
z,w\in\D,
$$
where $H_1(z,w)$ is the harmonic function of $z$ which is given by
$$
H_1(z,w) = (1-|w|^2)\bigg\{\frac12\,(3-|w|^2)\,\frac{1-|zw|^2}{|1-z\bar
w|^2}
+(1-|w|^2)\,\real\left[\frac{z\bar w}{(1-z\bar w)^2}\right]\bigg\},
$$
for $z,w\in\D$. Thus, $\widetilde\Gamma_1$ satisfies
$$
\Delta(1-|z|^2)^{-1}\Delta\widetilde\Gamma_1(z,w)=\delta_w(z),\qquad
z\in\D.
$$
We shall now establish the specified inequality for 
$\widetilde\Gamma_1$, which
entails that
$$\widetilde\Gamma_1(z,w)=O\big((1-|z|)^3\big),\qquad |z|\to1^-.$$
It follows that for a fixed $w\in\D$, the function $\widetilde\Gamma_1$
satisfies the boundary conditions
$$
\left\{\begin{array}{l}
\widetilde\Gamma_1(z,w)=0, \qquad z\in \T, \\
\partial_{n(z)}\,\widetilde\Gamma_1(z,w)=0, \qquad z\in \T,
\end{array}\right.
$$
where $\partial_{n(z)}$ is the outer normal derivative. In view of
this,
we conclude that $\widetilde\Gamma_1=\Gamma_1$. It remains to establish
the
claimed bounds for $\widetilde\Gamma_1$, from above and from below.
We use the following estimate of the logarithm:
$$\frac{r}2-\frac1{2r}<\log r<-\frac32+2r-\frac{r^2}2,\qquad 0<r<1,$$
with the choice
$$r=\left|\frac{z-w}{1-z\bar w}\right|^2,$$
and base our calculations on the identity
$$1-\left|\frac{z-w}{1-z\bar w}\right|^2=\frac{(1-|z|^2)(1-|w|^2)}
{|1-z\bar w|^2}.$$
The steps are rather lengthy but quite elementary, and are therefore
left to the reader as an exercise.
\end{proof}

The next result is a consequence of the positivity of the Green
function $\Gamma_1$.

\begin{lemma}
\label{lema intro}
Let $u$ be a $C^2$-smooth subharmonic function on $\D$ such that for
some real $\beta$, $0<\beta<2$,
$$|u(z)|=O\left(\frac1{(1-|z|)^\beta}\right)
\qquad\text{as}\quad |z|\to1^-.$$
Also, let $\varphi_\lambda=\varphi_{\lambda,1}$ denote the extremal
function for the point $\lambda\in\D$ in the space $\A_1^2(\D)$. Then
we have the following inequality:
$$
\int_\D u(z)(1-|z|^2)\,d\Sigma(z)\le
\int_\D|\varphi_\lambda(z)|^2u(z)(1-|z|^2)\,d\Sigma(z).
$$
\end{lemma}

\begin{proof}
We consider the potential function
$$
\Phi_\lambda(z)=\int_\D
G(z,w)\,\big(|\varphi_\lambda(w)|^2-1\big)(1-|w|^2)\,d\Sigma(w), \qquad
z\in\D,
$$
which solves the problem
$$
\left\{\begin{array}{l}
\Delta\Phi_\lambda(z)=\big(|\varphi_\lambda(z)|^2-1\big)(1-|z|^2),
\quad z\in \D,
\\[0.1cm]
\Phi_\lambda(z)=0,\quad z\in\T. \\[0.1cm]
\end{array}\right.
$$
Note that by using Green's formula as in \cite{Hed1}, we see that the
property that $\varphi_\lambda$ has according to Lemma \ref{lemma 2.1}
may be rephrased in terms of the potential function $\Phi_\lambda$:
$$\frac\partial{\partial n(z)}\Phi_\lambda(z)=0,\qquad z\in\T.$$
It follows that $\Phi_\lambda$ solves the over-determined boundary
value
problem
$$
\left\{\begin{array}{l} \Delta\Phi_\lambda(z) =
\big(|\varphi_\lambda(z)|^2-1\big)(1-|z|^2), \qquad z\in \D,
\\[0.3cm]
\Phi_\lambda(z)=0,\qquad z\in\T,
\\[0.1cm]
\dis\frac\partial{\partial n(z)}\Phi_\lambda(z)=0,\qquad z\in\T.
\end{array}\right.
$$
We may remove this over-determination by increasing the degree of the
elliptic operator:
$$
\left\{\begin{array}{l}
\dis\Delta\frac1{1-|z|^2}\Delta\Phi_\lambda(z) =
|\varphi_\lambda'(z)|^2, \qquad z\in \D,
\\[0.3cm]
\Phi_\lambda(z)=0,\qquad z\in\T,
\\[0.1cm]
\dis\frac\partial{\partial n(z)}\Phi_\lambda(z)=0,\qquad z\in\T.
\end{array}\right.
$$
This problem has a unique solution, which may be expressed in terms of
the Green function $\Gamma_1$,
$$
\Phi_\lambda(z)=\int_\D\Gamma_1(z,w)\,|\varphi'_\lambda(w)|^2d\Sigma(w)\geq0,
\quad z\in\D.
$$
We first apply this to the case when $u$ is $C^2$-smooth on $\bar\D$,
and see that in view of the assumption that $u$ is subharmonic, we
find,
by an application of Green's theorem, that
$$
\int_\D\big(|\varphi_\lambda(z)|^2-1\big)\,u(z)\,(1-|z|^2)\,d\Sigma(z)
=
\int_\D \Phi_\lambda(z)\,\Delta u(z)\, d\Sigma(z) \geq 0.
$$
If $u$ is not smooth up to the boundary, we perform the above for the
dilated function $u_r(z)=u(rz)$, with $0<r<1$. By the growth assumption
on
$u$ and the smoothness of the function $\varphi$ up to the boundary,
we may let $r\to1^-$ and apply Lebesgue's dominated convergence
theorem,
to conclude that
$$
\int_\D\big(|\varphi_\lambda(z)|^2-1\big)\,u(z)\,(1-|z|^2)\,d\Sigma(z)\geq
0.
$$
The proof is complete.
\end{proof}

Now, let $\omega$ be a strictly positive $C^\infty$-smooth weight on
$\bar\D$,
which is such that the function
$$z\mapsto \log\frac{\omega(z)}{1-|z|^2}$$
is subharmonic on $\D$. As before, let $K_\omega(z,w)$ denote the
reproducing
kernel for the weighted Bergman space $\A^2(\D,\omega)$. The following
result is basic for our further considerations.

\begin{thm}
The function $K_\omega$ extends to be $C^\infty$-smooth on the set
$(\bar\D\times\D)\cup(\D\times\bar\D)$.
\label{thm-ell-reg}
\end{thm}

This is a classical theorem of elliptic regularity type, obtained
in 1955 by L. Nirenberg \cite{Niren}. It is independent of
the above subharmonicity requirement.
\medskip

Next, we consider the function $\Lambda_\omega:\D\to\re$, as given by
$$
\Lambda_\omega(z)=\frac{|K_\omega(z,0)|^2}{K_\omega(0,0)}\,
\frac{\omega(z)}{2\,(1-|z|^2)}, \qquad z\in\D;
$$
it is positive and subharmonic in $\D$, due to the assumptions on
$\omega$,
and it has the growth behavior
$$\Lambda_\omega(z)=O\left(\frac1{1-|z|}\right),\qquad |z|\to1^-.$$
A basic property of $\Lambda_\omega$ is the following.

\begin{lemma}\label{lemma 2.1.1}
For each bounded harmonic function $h$ on $\D$, we have that
$$
\int_\D h(z)\,\Lambda_\omega(z)\,2\big(1-|z|^2\big)\,d\Sigma(z) = h(0).
$$
\end{lemma}

\begin{proof}
First, let us assume that $h$ is an analytic polynomial. Then, using the
reproducing property of the kernel function $K_\omega$, we have that
\begin{multline}
\label{eq-lemma2.1.1}
\int_\D h(z)\Lambda_\omega(z)\,2\,(1-|z|^2)\, d\Sigma(z)
\\
=\int_\D h(z)\frac{K_\omega(z,0)}{K_\omega(0,0)}\cj K_\omega(z,0)
\omega(z)\,d\Sigma(z)=h(0).
\end{multline}
Taking complex conjugates in (\ref{eq-lemma2.1.1}), we see that
the desired equality holds for all harmonic polynomials. An
approximation argument finishes the proof.
\end{proof}

This means that the function $\Lambda_\omega$ has a lot in common
with the function $|\varphi_\lambda|^2$ for the parameter
$\alpha=1$, which suggests it may have an expansive multiplier
property that is similar to the one obtained in %
Lemma~\ref{lema intro}.

\begin{lemma}
\label{lema expansivo}
Suppose that $u$ is a $C^2$-smooth subharmonic function in $\D$, which
has the growth bound
$$|u(z)|=O\left(\frac1{(1-|z|)^\beta}\right)
\qquad\text{as}\quad |z|\to1^-.$$
for some real $\beta$, $0<\beta<1$. We then have
$$
\int_\D u(z)\,2\,(1-|z|^2)\,d\Sigma(z)\le
\int_\D\Lambda_\omega(z)u(z)\,2\,(1-|z|^2)\,d\Sigma(z).
$$
\end{lemma}
\begin{proof}
We introduce the potential function
\begin{multline*}
\Phi_\omega(z)=\int_\D G(z,w)\,\big[\Lambda_\omega(w)-1\big]
\,2\,(1-|w|^2)\,d\Sigma(w)\\
=\int_\D G(z,w)\,
\bigg(\frac{|K_\omega(w,0)|^2}{K_\omega(0,0)}\,\omega(w)-2\,(1-|w|^2)\bigg)
\,d\Sigma(w),\qquad z\in\D,
\end{multline*}
which solves the boundary value problem
$$
\left\{\begin{array}{l}
\Delta\Phi_\omega(z)=2\,(1-|z|^2)\,\big[\Lambda_\omega(z)-1\big],
\qquad z\in\D, \\[0.1cm]
\Phi_\omega(z)=0, \qquad z\in\T.
\end{array}\right.
$$
It follows from Green's formula and Lemma \ref{lemma 2.1.1}, as in
\cite{Hed1}, that the potential function has
$$
\frac{\partial}{\partial n(z)}\Phi_\omega(z)=0,\qquad z\in\T.
$$
Then $\Phi_\omega$ solves the over-determinated boundary value problem
$$
\left\{\begin{array}{l}
\Delta\Phi_\omega(z)=2\,(1-|z|^2)\,\big[\Lambda_\omega(z)-1\big],
\qquad z\in\D
\\[0.1cm]
\Phi_\omega(z)=0, \qquad z\in\T,
\\[0.1cm]
\dis\frac\partial{\partial n(z)}\Phi_\omega(z)=0, \qquad z\in\T.
\end{array}\right.
$$
Increasing the degree of the elliptic operator, we find that
$\Phi_\omega$ also solves the problem
$$
\left\{\begin{array}{l} %
\Delta(1-|z|^2)^{-1}\Delta\Phi_\omega(z)=2
\Delta\Lambda_\omega(z)\geq0, \qquad z\in\D
\\[0.1cm]
\Phi_\omega(z)=0, \qquad z\in\T,
\\[0.1cm]
\dis\frac\partial{\partial n(z)}\Phi_\omega(z)=0, \qquad z\in\T,
\end{array}\right.
$$
which has a unique solution. Then $\Phi_\omega$ may be expressed in
terms of the Green function $\Gamma_1$,
$$
\Phi_\omega(z) =
2\int_\D\Gamma_1(z,w)\Delta\Lambda_\omega(w)\,d\Sigma(w)\geq0,
\qquad z\in\D.
$$
Let us consider now that $u$ is a subharmonic function which is
$C^2$-smooth on $\cj\D$. Then, by applying Green's theorem, it
follows that
$$
\int_\D(\Lambda_\omega(z)-1)\, u(z)\,2\,(1-|z|^2)\, d\Sigma(z) =
\int_\D\Phi_\omega(z)\Delta u(z)\, d\Sigma(z) \geq 0.
$$
In the case when $u$ is not smooth up to the boundary, we consider
the dilated function $u_r(z)=u(rz)$, with $0<r<1$, for which the
above inequality holds. It follows from the growth bounds of $u$
and $\Lambda$, and Lebesgue´s dominated convergence theorem, that
we may let $r\to 1^-$ to conclude that
$$
\int_\D(\Lambda_\omega(z)-1)\, u(z)\,2\,(1-|z|^2)\, d\Sigma(z) \geq 0.
$$
The proof is complete.
\end{proof}

\section{The proof of Theorem \ref{thm-main}}

We realized back in the introduction that, in order to obtain
Theorem \ref{thm-main}, all we need to do is show that
$$K_\omega(z,w)\neq0,\qquad (z,w)\in\big(\bar\D\times\D\big)\cup
\big(\D\times\bar\D\big),$$
provided that $\omega$ is a $C^\infty$-smooth and positive weight
function
on $\bar\D$, with the property that
$$z\mapsto\log\frac{\omega(z)}{1-|z|^2}$$
is subharmonic on $\D$. After all, Theorem \ref{thm-ell-reg}
guarantees that the weighted Bergman kernel is $C^\infty$-smooth
up to the boundary. It is easy to verify that the above
subharmonicity requirement is the same as the condition on the
curvature form in the statement of Theorem \ref{thm-main}.
\medskip

The proof splits naturally into two parts.

\begin{prop}\label{inside}
Under the above conditions on $\omega$,
$$K_\omega(z,w)\neq0,\qquad (z,w)\in\D\times\D.$$
\end{prop}

\begin{proof}
We observe first that for any M\oe{}bius map $\phi$ preserving the disk $\D$,
we have
$$
K_\omega(\phi(z),\phi(w))=K_{\omega_\phi}(z,w), \qquad z,w\in\D,
$$
where
$$
\omega_\phi(z)=|\phi'(z)|^2\, \omega\circ\phi(z),\qquad z\in\D.
$$
We claim that $\omega_\phi$ is a weight of the same type as
$\omega$. In fact, the function
$$
z\mapsto \log\bigg(\frac{\omega_{\phi}(z)}{1-|z|^2}\bigg)
$$
is subharmonic on $\D$ if, and only if, the function
$$
z\mapsto \log\bigg(\frac{\omega(z)}{1-|\phi^{-1}(z)|^2}\bigg)
$$
is subharmonic on $\D$ as well. Then, if we consider a M\oe{}bius map
$$
\phi^{-1}(z)=\gamma\frac{z-\zeta}{1-z\cj \zeta}, \qquad
z,\zeta\in\D, \quad |\gamma|=1,
$$
we find that
$$
\log\bigg(\frac{\omega(z)}{1-|\phi^{-1}(z)|^2}\bigg) =
\log\bigg(\frac{\omega(z)}{1-|z|^2}\bigg) +
\log\bigg(\frac{|1-z\cj\zeta|^2}{1-|\zeta|^2}\bigg).
$$
Thus, the function
$$
z\mapsto \log\bigg(\frac{\omega_{\phi}(z)}{1-|z|^2}\bigg)
$$
is subharmonic on $\D$ if, and only if,
$$
z\mapsto \log\bigg(\frac{\omega(z)}{1-|z|^2}\bigg)
$$
is subharmonic on $\D$ as well. It follows that it is enough to
specialize to $w=0$:
$$K_\omega(z,0)\neq0,\qquad z\in\D.$$
We note that, by the reproducing property of the kernel function,
$K_\omega(0,0)=\|K_\omega(z,0)\|_\omega^2>0.$
We introduce the extremal function $L$, given by
$$L(z)=(K_\omega(0,0))^{-\frac12}K_\omega(z,0), \qquad z\in\D,$$
which solves the problem
$$\sup\{\real f(0):\:\|f\|_\omega\leq1\}.$$
By Theorem~\ref{thm-ell-reg}, the function $L$ is
$C^\infty$-smooth on $\cj\D$.

We argue by contradiction. So, we assume that there exists a
$\lambda\in\D$ such that $K_\omega(\lambda,0)=0$; then
$L(\lambda)=0$. Consider then the function
$$\tilde L(z)=L(z)/\varphi_\lambda(z), \qquad z\in\D$$
where $\varphi_\lambda$ is the canonical divisor of $\{\lambda\}$
in the space $\A_1^2(\D)$. We should point out that, due to the
smoothness of $\varphi_\lambda$ and the fact that it doesn't have
any extraneous zeros on $\cj\D$, the function $\tilde L$ is also
$C^\infty$-smooth on $\cj\D$.

Let $u$ be the function given by
$$
u(z)=\frac{\omega(z)}{2(1-|z|^2)}|\tilde L(z)|^2, \qquad
z\in\cj\D.
$$
It follows from the hypothesis on the weight that $\log u$ is a
subharmonic function in $\D$, so that in particular, $u$ is
subharmonic as well and it has the growth bound
$$
|u(z)|= O\left(\frac1{1-|z|}\right), \qquad\text{as}\quad |z|\to1^-.
$$
It follows from Lemma \ref{lema intro} that
\begin{multline*}
\|\tilde L\|_\omega^2= \int_\D |\tilde L(z)|^2\,\omega(z)\, d\Sigma(z)
\\
=\int_\D u(z)\,2\,(1-|z|^2)\, d\Sigma(z)
\leq \int_\D|\varphi_\lambda(z)|^2\,u(z)\, 2\,(1-|z|^2)\, d\Sigma(z)
\\
= \int_\D |L(z)|^2\,\omega(z)\, d\Sigma(z) = \|L\|^2_\omega=1.
\end{multline*}
On the other hand, it follows from equation~(\ref{extremal
function}) that $\varphi_\lambda(0)<1$ and so $\tilde L(0)>L(0)$,
which violates the extremal property of $L(z)$. This is the
desired contradiction. Hence, the function $L$ does not have
zeroes in $\D$.
\end{proof}

\begin{prop}
Under the above conditions on $\omega$,
$$K_\omega(z,w)\neq0,\qquad (z,w)\in\big(\T\times\D\big)\cup
\big(\D\times\T\big).$$
\end{prop}

\begin{proof}
Note that, by the same argument used in the proof of the
Proposition~\ref{inside} and the fact that
$$
K_\omega(z,w)=\cj K_\omega(w,z), \qquad z,w\in\cj\D,
$$
it is enough to prove that
$$
K_\omega(z,0)\neq 0, \qquad z\in\T.
$$
We argue by contradiction. So, we shall assume that there exists
$\lambda\in\T$ such that $K_\omega(\lambda,0)=0$. Then, due to the
smoothness of $z\mapsto K_\omega(z,0)$ on $\bar \D$, as provided
by Theorem~\ref{thm-ell-reg}, we find that
$$
|K_\omega(z,0)|=O\big(|z-\lambda|\big), \qquad \text{as}\quad
\D\ni z \to\lambda.
$$
Let $0\leq r<1$ and define the function $f_r\in\A_1^2(\D)$, given by
$$
f_r(z)=\frac{K_1(z,r\lambda)}{\sqrt{K_1(r\lambda,r\lambda)}}=
\frac{(1-r^2)^{\frac32}}{(1-r\cj\lambda z)^3}, \quad z\in \D,
$$
where $K_1$ is the reproducing kernel for the space $\A_1^2(\D)$.
It follows that
$$
\|f_r\|_1=1, \qquad 0\leq r<1,
$$
where $\|\cdot\|_1$ is the norm in $\A_1^2(\D)$, as defined back
in Section 2. Furthermore, $|f_r|^2$ is bounded on $\D$ for each
$0\leq r<1$. We now consider the function
$$
R_\omega(z)=\frac{|K_\omega(z,0)|^2}{K_\omega(0,0)}\omega(z),
\quad z\in\D.
$$
It follows that
$$
R_\omega(z) = O\big( |z-\lambda|^2 \big)\qquad \text{as}\quad
\D\ni z \to\lambda.
$$
Then there exists a positive constant $M$ such that
$$
R_\omega(z)\leq M\, |z-\lambda|^2, \qquad z\in\D.
$$
It follows from Lemma~\ref{lema expansivo} and the inequality
$$
r\, |z-\lambda| \le \big|1-r\,\cj\lambda z\big|, \qquad z\in\D,
\quad 0<r<1,
$$
that
\begin{multline*}
1 = \int_\D|f_r(z)|^2\,2\,(1-|z|^2)\, d\Sigma(z) \leq
\int_\D\Lambda_\omega(z)|f_r(z)|^2\, 2\, (1-|z|^2)\, d\Sigma(z)
\\
=\int_\D R_\omega(z)\, |f_r(z)|^2\, d\Sigma(z)\leq M \int_\D
|z-\lambda|^2\, |f_r(z)|^2\, d\Sigma(z)
\\
\le M\,\frac{(1-r^2)^3}{r^2}\int_\D \frac1{|1-r\cj\lambda
z|^4}d\Sigma(z)= M\frac{1-r^2}{r^2},
\end{multline*}
which is a contradiction for $r$ sufficiently close to 1.
\end{proof}

\section{The proof of Theorem \ref{thm-nec}}

In this section, we construct, for a fixed $\alpha\geq\alpha_0=1.04$, 
an explicit example of a function $\mu\in C^\infty(\D)$, 
with the associated 2-form $\boldsymbol\mu=\mu d\Sigma$, such that
\begin{equation}
\boldsymbol\mu(z)+\frac\alpha2\, {\mathbf K}_{\mathbb
H}(z)\le0,\qquad z\in\D, \label{requerimiento}
\end{equation}
for which the optimization problem (OP) fails to have a
smooth-positive solution $\omega_0$ on $\D$, as the weighted
Bergman kernel $K_{\omega_1}(\cdot,0)$ has an extraneous
zero in $\D$. Here, $\omega_1$ is associated with $\mu$ as in the
introduction:
\begin{equation}
\log\omega_1(z)=-\int_\D \log\left|\frac{z-w}{1-z\bar w}\right|\,
\mu(w)\, d\Sigma(w), \qquad z\in\D.
\label{eq omega1}
\end{equation}
\medskip

\noindent{\bf The choice of $\boldsymbol\mu$.} We first consider the
extremal case for the inequality (\ref{requerimiento}),
$$
\boldsymbol\mu_{\mathbb H}(z) = -\frac\alpha2\, {\mathbf
K}_{\mathbb H}(z) = \frac{2\,\alpha\, d\Sigma(z)}{(1-|z|^2)^2},
\quad z\in\D.
$$
In this case we can compute explicitly the weight
$\omega_0$, which solves the optimization problem (OP),
$$
\omega_0(z)=(1-|z|^2)^\alpha, \qquad z\in\D.
$$
Following the same line of thought as in our toy example from back in the
introduction, we then consider the data function
\begin{equation}
\mu(z)=\frac{2\,\alpha}{(1-|z|^2)^2}-\sum_{k}
\rho_k\,\delta_{a_k}(z), \qquad z\in\D, \label{eq mu}
\end{equation}
where $A=\{a_k\}_k$ is a finite collection of points in $\D, \;
\{\rho_k\}_k$ is a convenient sequence of positive constants, and
$\delta_{a_k}$ stands for the Dirac delta function, concentrated
at the point $a_k\in\D$. The 2-form $\boldsymbol\mu=\mu d\Sigma$
meets the inequality (\ref{requerimiento}), but $\mu$ is very rough 
at the points of $A$ and hence it does not satisfy the $C^\infty$-smoothness
requirement. However, it is easy to approximate a point-mass by a sequence
of positive $C^\infty$-smooth functions. Moreover, we may replace the
function 
$$\nu(z)=\frac{2\,\alpha}{(1-|z|^2)^2}$$
by a slight dilation,
$$\nu_r(z)=\frac{2\, r^2\,\alpha}{(1-r^2\,|z|^2)^2}$$
with $r$, $0<r<1$, close to $1$. This means that if 
\begin{equation}
K_{\omega_1}(z,w)\neq 0, \qquad z,w\in\D,
\label{K no zero}
\end{equation}
holds for smooth indata $\mu$, with $\mu$ and $\omega_1$ connected via 
(\ref{eq omega1}), then it also holds for rough indata of
the above type, modulo some slight modifications. To explain these
modifications, we note that 
$$
\int_\D \log\left|\frac{z-w}{1-z\bar w}\right|\,\nu_r(w)\, d\Sigma(w) =
\alpha\log\frac{1-r^2}{1-r^2|z|^2}, \qquad z\in\D.
$$
Let the weight $\omega_{1,r}$ be defined by the formula (\ref{eq omega1}),
where $\omega_1$ is replaced
by $\omega_{1,r}$, and $\mu$ is replaced by
$$
\mu_r(z)=\frac{2\,r^2\,\alpha}{(1-r^2\,|z|^2)^2}-\sum_{k}
\rho_k\,\delta_{a_k}(z), \qquad z\in\D.
$$
We then calculate that
$$
\omega_{1,r}(z)=(1-r^2)^{-\alpha}\big(1-r^2\,|z|^2\big)^\alpha
\prod_{k}\left|\frac{z-a_k}{1-\bar a_kz}\right|^{\rho_k},
\qquad z\in\D.
$$
As reproducing kernel functions have homogeneity index $-1$ in general,
that is,
$$
K_{t\omega}(z,w)=\frac1tK_\omega(z,w),
$$
holds for arbitrary positive constant $t$, we see that
$$
K_{\omega_{1,r}}(z,w)=(1-r^2)^\alpha\, K_{\omega_{2,r}}(z,w), 
$$
where $\omega_{2,r}$ is the weight
$$
\omega_{2,r}(z)=\big(1-r^2\,|z|^2\big)^\alpha
\prod_{k}\left|\frac{z-a_k}{1-\bar a_kz}\right|^{\rho_k},
\qquad z\in\D.
$$
If (\ref{K no zero}) holds for $\omega_1=\omega_{1,r}$, then it also holds
for $\omega_1=\omega_{2,r}$, and vice versa. As $r\to 1^-$, the weight
$\omega_{2,r}$ tends to
\begin{equation}
\omega_{2}(z)=\big(1-|z|^2\big)^\alpha
\prod_{k}\left|\frac{z-a_k}{1-\bar a_kz}\right|^{\rho_k},
\qquad z\in\D,
\label{eq omega2}
\end{equation}
and if the reproducing kernel functions for the weights $\omega_{2,r}$
are all zero-free in $\D^2$, then so is the reproducing kernel for the weight
$\omega_2$, by a limit process argument. We shall prove that with appropriate 
choices of the configuration 
of the points $A=\{a_k\}_k$ as well as of the positive parameters $\rho_k$, 
the reproducing kernel function for $\omega_2$ will have zeros in $\D^2$.
This then shows that also with smooth data, we must have zeros in the 
associated reproducing kernel function.

We should note that $\omega_2$ satisfies
$$
\Delta\log\omega_2(z)=-\frac12\,\mu(z),\qquad z\in\D,
$$
where $\mu$ is given by (\ref{eq mu}). This means that $\omega_2$ is the
weight we where looking for.

For computational reasons, we shall only consider $\rho_k$ such 
that $\rho_k/2$ is a positive integer.

\medskip

\noindent{\bf Relations between kernel functions.} A closed subspace 
$I$ of $\A_\alpha^2(\D)$ is called
invariant if it is invariant under the multiplication by the
identity function, more precisely, if $zf\in I$ whenever $f\in I$.
For a sequence $A=\{a_1,\,a_2,\,\ldots\}$ of points on $\D$, the
subspace $I_A$ consisting of all functions in $A^2_\alpha(\D)$
whose zero sets contain $A$, counting multiplicities, is an
invariant subspace. Such a subspace $I_A$ is called {\sl zero-based
invariant subspace}. The reproducing kernel function for the 
invariant subspace $I_A$ is as usual defined by the formula
$$
K^\alpha_A(z,w)=\sum_{n=0}^{+\infty} e_n(z)\,\bar e_n(w), \qquad z,w\in\D,
$$
where the functions $e_1(z),\, e_2(z),\, e_3(z), \ldots$ form an orthonormal
basis for $I_A$. It has the reproducing property
$$
f(z)=\int_\D K^\alpha_A(z,w)\, f(w) \, (\alpha+1)\big(1-|w|^2\big)^\alpha 
d\Sigma(w),\qquad z\in\D,
$$
for $f\in I_A$. For a finite subset $A=\{a_k\}_k$ of $\D$, the associated
(finite) Blaschke product is the function
$$
B_A(z)=\prod_{k}\frac{z-a_k}{1-\bar a_k z}, \qquad z\in\D.
$$
We consider the following two weights: 
$$
\omega_\alpha(z)=(\alpha+1)\big(1-|z|^2\big)^\alpha, \qquad
\omega_{\alpha,A}(z)=(\alpha+1)\big(1-|z|^2\big)^\alpha\, |B_A(z)|^2. 
$$
Note that $\omega_{\alpha,A}$ equals the weight $\omega_2$ as defined by 
(\ref{eq omega2}), with all the parameters $\rho_k$ set equal to $2$.

The following proposition is well known.

\begin{prop}
We have the following identity of kernels:
$$
K^\alpha_A(z,w)=B_A(z)\,\bar B_A(w)\, K_{\omega_{\alpha,A}}(z,w), \qquad
z,w\in\D.
$$
\end{prop}

In view of the above proposition, we need to look for extraneous zeros 
in the reproducing kernel function for $I_A$ in order to get zeros of
the kernel function for the weight $\omega_2=\omega_{\alpha,A}$.

\medskip
\noindent{\bf The kernel function for a zero-based invariant
subspace.} 
 When the sequence $A$ consist of a finite
number of distinct points, the kernel function $K_A^\alpha$ for
$I_A$ may be obtained by means of the well-known iterative formula
(see \cite{HHopen})
\begin{eqnarray}
\dis K_\emptyset^\alpha(z,w)=\frac 1{(1-z\cj w)^{(\alpha+2)}},
\qquad z,w \in\D
\label{Kvacio}
\\[.25cm]
\dis K^\alpha_{A\cup\{\lambda\}}(z,w) = K^\alpha_{A}(z,w) -
\frac{K^\alpha_{A}(z,\lambda)K^\alpha_{A}(\lambda,w)}{K^\alpha_{A}(\lambda,\lambda)},
\qquad \lambda\notin A. \label{KA+b}
\end{eqnarray}
The first step of this iteration is to apply the formula~(\ref{KA+b})
to the case $A=\emptyset$ and $\lambda\in\D$, to get
$$
K_\lambda^\alpha(z,w)=\frac1{(1-z\bar w)^{\alpha+2}}-
\frac{(1-|\lambda|^2)^{\alpha+2}}{(1-z\bar\lambda)^{\alpha+2}(1-\lambda\bar
w)^{\alpha+2}}, \qquad z,w \in\D.
$$
In the case where the zero set $A$ contains repeated points, i.e.
zeros with high multiplicity, the kernel function $K^\alpha_A$ for
$I_A$ can be computed via an iterative formula, similar to
(\ref{KA+b}), but involving the derivatives of the kernel. Namely, if we
assume $A$ has no multiple points, then
\begin{eqnarray*}
\dis K_\emptyset^\alpha(z,w)=\frac 1{(1-z\cj w)^{(\alpha+2)}},
\qquad z,w \in\D,
\\[.25cm]
\dis K^\alpha_{A\cup\{\lambda\}}(z,w) = K^\alpha_{A}(z,w) -
\frac{K^\alpha_{A}(z,\lambda)K^\alpha_{A}(\lambda,w)}{K^\alpha_{A}(\lambda,\lambda)},
\qquad \lambda\notin A,
\\[.25cm]
\dis K^\alpha_{A\cup\{\xi\}}(z,w)=K^\alpha_{A}(z,w) - \frac{\left.
\partial_z K^\alpha_{A}(z,w)\right|_{z=\xi}\, \left. \bar\partial_w
K^\alpha_{A}(z,w)\right|_{w=\xi}}{\left.\partial_z\bar\partial_w
K^\alpha_{A}(z,w)\right|_{z=\xi,\: w=\xi}}, \qquad \xi\in A,
\end{eqnarray*}
where
$$
\partial_z=\frac12\,\left(\frac\partial{\partial x} -\ii \frac\partial{\partial
y}\right), \qquad
\bar\partial_z=\frac12\,\left(\frac\partial{\partial x} +\ii
\frac\partial{\partial y}\right), \qquad z=x+\ii y.
$$
When $A$ consist of $n$ copies of the point $a\in\D$,
the reproducing kernel function $K^\alpha_A$ for the subspace $I_A$ is
\begin{multline*}
K^\alpha_A(z,w)=\frac1{(1-z\bar w)^{\alpha+2}} \\
-\frac{(1-|a|^2)^{\alpha+2}}{(1-z\bar a)^{\alpha+2}(1-a\bar
w)^{\alpha+2}}\, \sum_{j=0}^{n-1}\, C_{\alpha,j}
\,
\Big(\frac{z-a}{1-z\bar a}\Big)^j\Big(\frac{\bar w -\bar
a}{1-a\bar w}\Big)^j,
\end{multline*}
where 
$$
C_{\alpha,j}=\dis \frac{\Gamma(\alpha+2+j)}{\Gamma(\alpha+2)j!}.
$$

Our computational work suggest that whenever $\alpha$ is bigger than $3$,
it is possible to find extraneous zeros; for $\alpha\leq 3$, however,
it seems that such extraneous zeros do not occur.
\medskip

\noindent{\bf Computational implementation.} The recursive formula for 
(\ref{KA+b}) is somewhat inconvenient in computational applications. Instead,
we use the fact that the kernel function with $w=0$ may be expressed as
\begin{equation}
K_A^\alpha(z,0)=1-\sum_{j}\frac{c_j}{(1-z\cj
a_j)^{\alpha+2}}, \label{kernelconst}
\end{equation}
where the $\{c_j\}_j$ are certain complex coefficients, which depend on 
the sequence$A$ and on the parameter $\alpha$. To determine them, we need 
to solve the linear system $K_A^\alpha(a_j,0)=0$, for every $a_j\in A$. This 
we rewrite in terms of a square matrix $\mathbf M$, whose entries
are given by
$$
M_{i,j}=\frac1{(1-a_i\cj a_j)^{\alpha+2}},\qquad
a_i,a_j\in A;
$$
the equation that determines the coefficients is the linear system
\begin{equation}
{\mathbf M}\, \mathbf{c} = \mathbf{1}, %
\label{sistema lineal}
\end{equation}
where $\mathbf{c}=\{c_j\}_j$ is the coefficient vector in column form, 
and $\mathbf{1}$ is the column vector with the number
$1$ in all its coordinates.

Although the numerical computation of the kernel function
$K_A^\alpha$ for several zeros via the linear system (\ref{sistema lineal}) 
is technically possible, we are to treat a rather badly scaled
matrix. Numerically, the matrix $\mathbf M$ is
nearly singular, and this leads to possibly large numerical error while
solving the system~(\ref{sistema lineal}). In addition, the
kernel function expressed as expressed by (\ref{kernelconst}) is quite
sensitive to even small perturbations in the appearing coefficients already 
when
the number of zeros in $A$ is rather modest. Taken together, this forces us 
to work with higher precision than what is standard. We use 40 decimal digits 
of precision, instead of the standard 16. Most of the numerical tests were 
done in \texttt{MATLAB} 6.1. To work with forty decimal digits, we used the 
command \texttt{VPA} (variable precision arithmetic).
We estimate the computational error in the approximated kernel function 
$K_A^\alpha(\cdot,0)$ by comparing the computed values at the points of 
the given zero set $A=\{a_k\}_k$, where the function vanishes. If these 
values are sufficiently small, we may be certain that the negativity 
of $K_A^\alpha(1,0)$ is a genuine phenomenon, provided that the negative 
value is substantially bigger than the numerically obtained values at the 
given collection of zeros.  
\medskip

\noindent {\bf The appearance of an extraneous zero.} We first observe 
that if the elements of the zero set $A$ are distributed symmetrically 
respect to the real axis, then the kernel function $K_A^\alpha(x,0)$ is 
real-valued for real $x$. Since, trivially, $0<K_A^\alpha(0,0)$, and since 
the function $K_A^\alpha(z,0)$ is continuous in the closed unit disk
$\bar\D$, it follows from the Mean Value Theorem of Calculus that 
$K_A^\alpha(x_0,0)=0$ holds for some $x_0$, $0<x_0<1$, provided that
$K_A^\alpha(1,0)<0$.

After testing several patterns for the distribution of the zero
set $A$ (including, for instance, a multiple zero at a single point), 
we focused our computations on a configuration which yields 
extraneous zeros for $\alpha$ all the way down to 1.04. The same pattern
seemed to emerge also when the computer was allowed to pick the configuration
of the given zeros according to a so-called {\sl genetic algorithm}. 
The configuration depends on the number $n$ of points of $A$, as well as on
two parameters $\theta$ and $d$, with $0<\theta<\frac12\pi$ and $1<d<+\infty$.
By the construction, $n$ is an even number. The points of $A$ are given by
$$
\begin{array}{l}
a_k=\exp \big\{3(\ii-\theta)\,d^{k-n/2}\big\},\qquad k = 1, 2, \ldots,
n/2, \\[0.2cm]
a_{k+\frac n2}=\bar a_k, \qquad\qquad\qquad\qquad k = 1, 2, \ldots, n/2.
\end{array}
$$

In Figure~\ref{fig1}, we plot the configuration of the given zero set in
the complex plane for $n=80$ zeros, with the parameter values $\alpha=1.5$,
$\theta=0.15$, and $d=1.3$.

\begin{figure}
\begin{center}
\caption{Zero set $A$ for $n=80$ generated by the parameters
$\theta=0.15, \; d=1.30$, for $\alpha=1.50$. The indicated straight 
lines are tangents to the curves on which the zeros are located.} 
\label{fig1} %
\vspace{-0.5cm}
\end{center}
\end{figure}

In Figure~\ref{fig2} (top), we plot the level curves around the extraneous 
zero of the modulus of the kernel function for the zero set $A$ generated 
by the parameter values  $\alpha=3$, $n=6$, $\theta=0.51$, and $d=10$.
In Figure \ref{fig2} (bottom), we plot the kernel function $K_A^\alpha(1,0)$, 
while varying $\alpha$ in the interval $[1.035,1.050]$.

Inspiration for our numerical work was derived from \cite{SJak}, where 
Jakobsson obtained extraneous zeros for $\alpha\approx 1.40$.

\begin{figure*}
\begin{center}
\caption{(top) Level curves of the kernel function $K_A^\alpha(\cdot,0)$
for the zero set $A$ generated by the parameter values: $\alpha=3,\;
n=6,\; \theta=0.51,\;d=10$.
(bottom) Values of the kernel function $K_A^\alpha(1,0)$ for the
zero set $A$ generated by the parameters: $n=1500,\;
\theta=0.033,\;d=1.045$ and $1.035 \leq \alpha\leq 1.050$} \label{fig2}%
\end{center}
\end{figure*}
\vfill\eject

In the table below (Table 1), we supply values of $n$, $\theta$, and $d$, 
for which our computations indicated that the kernel function possesses 
extraneous zeros for a prescribed value of the parameter $\alpha$.

\begin{center}
\begin{tabular}{|c c c c|}\hline
 $\alpha$ & $n$   & $\theta$  &  $d$       \\ \hline
    3     &  6    & 0.51    &  10      \\ \hline
   2.5    &  8    & 0.48    &  8       \\ \hline
    2     &  14   & 0.351   &  3       \\ \hline
   1.6    &  26   & 0.265   &  2.1     \\ \hline
   1.25   &  78   & 0.176   &  1.52    \\ \hline
   1.118  &  230  & 0.104   &  1.22    \\ \hline
   1.1072 &  272  & 0.092   &  1.183   \\ \hline
   1.097  &  340  & 0.07725 &  1.141   \\ \hline
   1.065  &  550  & 0.07    &  1.13    \\ \hline
   1.053  &  770  & 0.0556  &  1.09    \\ \hline
   1.046  &  944  & 0.0497  &  1.078   \\ \hline
   1.043  &  1090 & 0.0445  &  1.067   \\ \hline
   1.04   &  1500 & 0.033   &  1.045   \\ \hline
\end{tabular}
\\[0.5cm]
{\bf Table 1.}
\end{center}

\medskip

\noindent{\bf Hinted analytical solution of the problem.} The configuration
of zeros in the above numerically-based counterexamples which works at least
down to the parameter value $\alpha\approx1.04$ suggests that analytically,
we should ``smear out'' the zeros along the two lines; we may decide to work
almost infinitesimally close to to the point $1$, and by blowing up, we
may assume that the domain is the upper half plane, where the origin plays
the role of the point $1$. We place hyperbolically equi-distributed 
smeared-out ``zeros'' along two half-lines emanating from the origin, 
symmetrically located with respect to the imaginary axis. In geometric terms,
this means that we construct a new weight -- the old one being 
$(\,\text{Im}\,z)^\alpha$ -- which is the same as the old weight in the
region between the two half-lines, but is the old weight times the exponential
of a constant times the angle to the nearest half-line in the remaining
two regions near the real line. We suspect that the reproducing kernel
for this new weight (with one argument fixed equal to $\ii$) has a zero 
somewhere along the imaginary axis for each fixed $\alpha$, 
$1<\alpha<+\infty$, provided the angle of the two half-lines to the real 
line are chosen appropriately, and the constant that regulates the density
of the ``zeros'' is appropriate as well. However, the actual implementation
of this scheme is not easy, because it is generally a hard problem to 
calculate reproducing kernel functions for weights that do not exhibit strong
symmetry properties. We would like to be able to return to this problem.



\end{document}